\documentclass[a4paper, oneside, 11pt]{article}
\usepackage[english]{babel}
\usepackage{amsfonts}
\usepackage{amssymb}
\usepackage{graphics}
\usepackage{amsrefs}
\usepackage{latexsym}
\begin{document}

\title{{\bf{\Large{A new parametrization of the Gnedin-Fisher\\species sampling model}}}\footnote{{\it AMS (2000) subject classification}. Primary: 60G58. Secondary: 60G09.} }
\author{\textsc {Annalisa Cerquetti}\footnote{Corresponding author, SAPIENZA University of Rome, Via del Castro Laurenziano, 9, 00161 Rome, Italy. E-mail: {\tt annalisa.cerquetti@gmail.com}}\\
\it{\small Department of Methods and Models for Economics, Territory and Finance}\\
  \it{\small Sapienza University of Rome, Italy }}
\newtheorem{teo}{Theorem}
\date{\today}
\maketitle{}

\begin{abstract}

We introduce a new parametrization for the two-parameter species sampling model with {\it finite} but {\it random} number of different species recently introduced in Gnedin (2010a). We show the reparametrization yields a representation in terms of generalized Waring mixture of Fisher species sampling models and derive the structural distribution of the model.\\

\end{abstract}

\section{Introduction}

Gnedin (2010a) introduces a two parameter family of exchangeable partition models belonging to the  Gibbs class of genius $\alpha =-1$ (Gnedin and Pitman, 2006), by suitably mixing Fisher's (1943) $(-1, \xi)$ partitions over the fixed number of boxes. The resulting $(\gamma, \zeta)$ {\it Gnedin-Fisher} species sampling model has exchangeable partition probability function (EPPF) of the form
\begin{equation}
\label{gnedintwo}
p_{\gamma, \zeta}(n_1, \dots, n_k)=\frac{(\gamma)_{n-k} \prod_{i=1}^{k-1}(i^2 -\gamma i +\zeta)}{\prod_{l=1}^{n-1} (l^2 +\gamma l +\zeta)} \prod_{j=1}^k n_j!
\end{equation}
obtained by sequential construction with {\it one-step} allocation rules
$$
\mbox{(O)}: p_{j}({\bf n}) :=\frac{(n-k+\gamma)(n_j +1)}{n^2 +\gamma n +\zeta}      \mbox{ \hspace{0.5cm} for $j=1,\dots, k$  \hspace{0.5cm} and  \hspace{0.5cm} }     \ \mbox{(N)}: p_0({\bf n}):=\frac{k^2 -\gamma k +\zeta}{n^2 +\gamma n +\zeta},
$$
for $\gamma \geq 0$ and (i) either $i^2 -\gamma i +\zeta$ (strictly) positive for all $i \in \mathbb{N}$ or (ii) the quadratic is positive for $i \in \{1,\dots, i_0-1\}$ and has a root at $i_0$. (See Pitman, 2006, Hansen and Pitman, 2000, for background on exchangeable partitions and sequential constructions).  As from Theorem 1. in Gnedin (2010a) model (\ref{gnedintwo}) arises by mixing Poisson-Dirichlet $(-1, \xi)$ models (Pitman and Yor, 1997), for $\xi=1,2,3,\dots,$ 
over $\xi$ with 
\begin{equation}
\label{mix}
\mathbb{P}_{\gamma, \zeta} (\Xi=\xi)=\frac{\Gamma(z_1+1)\Gamma(z_2 +1)}{\Gamma(\gamma)}\frac{\prod_{i=1}^{\xi -1}(i^2 -\gamma i +\zeta)}{\xi!(\xi-1)!},
\end{equation}
for some complex $z_1$ and $z_2$. For $\zeta =0$ then $\gamma \in (0,1)$ (cfr. Gnedin, 2010a, Sect. 6) and (\ref{gnedintwo}) reduces to
\begin{equation}
\label{gnedin}
p_{\gamma}(n_1, \dots, n_k)= \frac{(k-1)!}{(n-1)!}\frac{(1-\gamma)_{k-1} (\gamma)_{n-k}}{(1+\gamma)_{n-1}} \prod_{j=1}^k n_j!.
\end{equation}
The sequence of the number of occupied boxes $K_n$ for both model (\ref{gnedin}) and (\ref{gnedintwo}) is a nondecreasing Markov chain with $0-1$ increments and transition probabilities determined by the specific rule (N), whose distribution follows by the general formula for Gibbs partitions of genius $\alpha \in (-\infty, 1)$
\begin{equation}
\label{nblocks}
\mathbb{P}(K_n=k)= V_{n,k} S_{n,k}^{-1, -\alpha},
\end{equation}
for $V_{n,k}$ the general Gibbs weights satisfying the backward recursion $V_{n,k}=(n-k\alpha)V_{n+1, k}+ V_{n+1, k+1}$ and $S_{n,k}^{-1, \alpha}$ generalized Stirling numbers. For $\alpha=-1$ those reduce to Lah numbers (see e.g. Charalambides, 2005) $S_{n,k}^{-1,1}= {n-1 \choose k-1} \frac{n!}{k!}$ hence, e.g. for the one-parameter model
\begin{equation}
\label{number}
\mathbb{P}_{\gamma}(K_n=k)= {n \choose k} \frac{(1-\gamma)_{k-1} (\gamma)_{n-k}}{(1+\gamma)_{n-1}}.
\end{equation}
As from Gnedin (2010a, cfr. eq. (9) and (10)) the mixing law yielding the one-parameter model (\ref{gnedin}) arises from (\ref{number}) for $n \rightarrow \infty$ by the standard asymptotics $\Gamma(n+\alpha)/\Gamma(n+b) \sim n^{a-b}$ and corresponds to  
\begin{equation}
\label{priorone}
\mathbb{P}_\gamma(\Xi=\xi)=\frac{\gamma(1-\gamma)_{\xi-1}}{\xi!},
\end{equation}
for $\xi=1,2,\dots,$ and $\gamma \in (0,1)$, while, for $1 \leq k \leq n$, a {\it posterior} distribution for $\Xi$ results 
\begin{equation}
\label{postgamma}
\mathbb{P}_\gamma(\Xi=\xi |K_n=k)=\frac{(n-1)!}{(k-1)!} \frac{\Gamma(\gamma +n)}{\Gamma(\gamma +n -k)} \frac{(k- \gamma)_{\xi} \Gamma(k+\xi)}{\Gamma(\xi +1)\Gamma(k+\xi+n)}.
\end{equation} 
\section{The new parametrization}
Gnedin (2010a, cfr. Sect. 2) points out that the Gibbs weights of model (\ref{gnedintwo}) 
can be split in linear factors by factoring the quadratics as
$$
x^2+\gamma x +\zeta =(x+z_1)(x+z_2), \mbox{  and  } x^2 -\gamma x + \zeta=(x+s_1)(x+ s_2) 
$$
thus providing the alternative {\it five} parameters representation 
\begin{equation}
\label{essezeta}
V_{n,k}^{\gamma, \zeta}= \frac{(\gamma)_{n-k}(s_1+1)_{k-1} (s_2 +1)_{k-1}}{(z_1 +1)_{n-1}(z_2 +1)_{n-1}}
\end{equation}
for some complex $z_1,z_2,s_1,s_2$, such that $z_1+z_2=\gamma$, $z_1z_2=\zeta$, $s_1 +s_2=-\gamma$, $s_1s_2=\zeta$. \\

Here we show those constraints limit admissible values for the four parameters $s_1, s_2, z_1$ and $z_2$ yielding an interesting alternative {\it two} parameter representation of the weights of the Gnedin-Fisher model.\\\\
{\bf Theorem 1.} {\it For $\psi \in [0,1)$ and $0 <\gamma < \psi +1$ the EPPF of the two-parameter $(\gamma, \zeta)$-Gnedin-Fisher species sampling model} (\ref{gnedintwo}) {\it admits the following alternative representation
\begin{equation}
\label{miogned}
p_{\gamma, \psi}(n_1, \dots, n_k)= \frac{(\gamma)_{n-k} (1- \psi)_{k-1} (1 -\gamma +\psi)_{k-1}}{(1+\psi)_{n-1}(1 +\gamma -\psi)_{n-1}} \prod_{j=1}^k n_j!.
\end{equation}
For $\psi=0$, then $\gamma \in (0,1)$  and} (\ref{miogned}) {\it yields the one-parameter Gnedin-Fisher model} (\ref{gnedin}).\\\\
{\it Proof:}
For $z_1+z_2=\gamma$, $z_1z_2=\zeta$, $s_1 +s_2=-\gamma$ and $s_1s_2=\zeta$ vectors $(z_1, z_2)$ and $(s_1, s_2)$ must be the roots (complex or real) of the following quadratic polynomials
$$
z_1^2 -\gamma z_1 +\zeta=0 \mbox{  \hspace{0.8cm}  and  \hspace{0.8cm}  }
z_2^2 -\gamma z_2 +\zeta=0,
$$
$$
s_1^2 +\gamma s_1 +\zeta=0  \mbox{ \hspace{0.8cm}    and \hspace{0.8cm}   }
s_2^2 +\gamma s_2+\zeta=0.
$$
For $\gamma^2 -4\zeta>0$ admissible {\it real} solutions are 
$$
z_1= \frac{\gamma \pm  \sqrt{\gamma^2 -4 \zeta}}{2}  \mbox{ \hspace{0.8cm}  and  \hspace{0.8cm} } z_2= \frac{\gamma \pm \sqrt{\gamma^2 -4 \zeta}}{2},
$$
$$
s_1= \frac{-\gamma \pm \sqrt{\gamma^2 -4\zeta}}{2} \mbox{\hspace{0.8cm}    and \hspace{0.8cm}} s_2 =\frac{ -\gamma \pm  \sqrt{\gamma^2 -4\zeta}}{2}.
$$
For $\gamma^2 -4\zeta<0$ admissible {\it complex} solutions are
$$
z_1= \frac{\gamma \pm i \sqrt{4 \zeta-\gamma^2}}{2}  \mbox { \hspace{0.8cm} and  \hspace{0.8cm} } z_2= \frac{\gamma \pm i \sqrt{4 \zeta-\gamma^2}}{2},
$$
$$
s_1= \frac{-\gamma \pm i \sqrt{4 \zeta-\gamma^2}}{2} \mbox{  \hspace{0.8cm}  and \hspace{0.8cm}} s_2 =\frac{ -\gamma \pm i \sqrt{4 \zeta-\gamma^2}}{2}.
$$
Now let indifferently $A=i\sqrt{4 \zeta -\gamma^2}/2$ or $A= \sqrt{\gamma^2 -4\zeta}/2$. Then, regardless of the solutions being real or complex, possible vectors satisfying the constraints $z_1+z_2=\gamma$, $s_1+s_2=-\gamma$, $z_1z_2= \zeta$ and $s_1s_2=\zeta$ must be as follows
$$
(z_1, z_2)=\left(\frac{\gamma}{2}+A, \frac{\gamma}{2}-A \right) \mbox{   or   } \left(\frac{\gamma}{2}-A, \frac{\gamma}{2}+A \right) 
$$
and 
$$
(s_1, s_2)=\left(-\frac{\gamma}{2}+A, -\frac{\gamma}{2}-A \right) \mbox{   or   } \left(-\frac{\gamma}{2}-A, -\frac{\gamma}{2}+A \right),
$$
which shows admissible solutions reduce to $z_1=-s_1$ and $z_2=-s_2$ or $z_1=-s_2$ and $z_2=-s_1$. Since (\ref{essezeta}) is invariant to permutations of $(z_1, z_2)$ and $(s_1, s_2)$ the five parameters in (\ref{essezeta}) reduce to $\psi$ and $\gamma$ for $z_1=\psi$ and $z_2= \gamma -\psi$, and $s_1=-\psi$, $s_2=\psi -\gamma$ thus yielding (\ref{miogned}). Moreover the positiveness of the numerator in (\ref{miogned}) implies $1-\psi>0$ and $1-\gamma +\psi>0$. For $\psi=0$, $0 < \gamma < 1$ and (\ref{miogned}) reduces to (\ref{gnedin}) by standard combinatorial calculus. \hspace{3.6cm}$\square$\\\\
{\bf Remark 2.} We stress that, as from Lemma 5.2 in Gnedin (2010b), beside the extended two-parameter Poisson-Dirichlet $(\alpha, \theta)$ family of partitions models, with $\alpha \in (0,1)$ and $\theta > -\alpha$ or $\alpha < 0$ and $\theta= |\alpha|\xi$, $\xi=1, 2,\dots$, the two-parameter Gnedin-Fisher models is the unique class of exchangeable random partitions with Gibbs weights in the nice multiplicative form 
\begin{equation}
\label{niceform}
V_{n,k}= \frac{\prod_{i=0}^{n-k-1}g_0(i) \prod_{j=1}^{k-1} g_1(j)}{\prod_{l=1}^{n-1} g(l)}
\end{equation}
for $g_0(\cdot)$, $g_1(\cdot)$ and $g(\cdot): \mathbb{N} \rightarrow \mathbb{R}$  satisfying the identity 
$$
(n- \alpha k)g_0(n-k)+g_1(k)= g(n), \hspace{0.5cm}  1 \leq k \leq n, n \in \mathbb{N}.
$$
In terms of equation (\ref{niceform})  the weights in (\ref{miogned}) may be written as 
$$
V_{n,k}^{\gamma, \psi}= \frac{\prod_{i=0}^{n-k-1}(i+\gamma) \prod_{j=1}^{k-1} (j-\psi)(j-\gamma +\psi)}{\prod_{l=1}^{n-1} (l +\psi)(l +\gamma -\psi)}
$$
for $g_0(i)=(i+\gamma)$, $g_1(j)=(j-\psi)(j -\gamma +\psi)$ and $g(l)=(l +\psi)(l +\gamma -\psi)$.\\\\
For the reparametrized model (\ref{miogned})  {\it multistep allocation rules} may be derived specializing the general form for Gibbs partitions of genius $\alpha \in (-\infty, 1)$  introduced in Cerquetti (2008) as follows. Start with box $B_{1,1}$, containing a single ball 1. At step $n$ the allocation of $n$ balls is a certain random partition $\Pi_n=(B_{n,1},\dots, B_{n,K_n})$ of the set of balls $[n]$. Given the number of boxes is $K_n=k$ and the occupancy counts are $(n_1,\dots, n_k)$, the partition of $[n+m]$ at step $n+m$ is  obtained by randomly placing the additional $m$ balls \\\\
(AO): in $k$ {\it old} boxes in configuration $(m_1,\dots, m_k)$, for $m_j \geq 0$, $\sum_{j=1}^k m_j=m$, with probability 
\begin{equation}
\label{gibbsallold}
p_{{\bf m}}({\bf n}):=\frac{(\gamma +n -k)_{m}}{(\psi +n)_m(\gamma -\psi +n)_m} \prod_{j=1}^k (n_j+1)_{m_j\uparrow},
\end{equation}
(AN): in $k^*$ {\it new} boxes in configuration $(s_1, \dots, s_{k^{*}})$, for $\sum_{j=1}^{k^*} s_j =m$, $1 \leq k^* \leq m$, $s_j \geq 1$, with probability  
\begin{equation}
\label{gibbsallnew}
p_{{\bf s}}({\bf n}):=\frac{(\gamma +n -k)_{m-k^*} (k- \psi)_{k^*} (k-\gamma +\psi)_{k^*}}{(\psi +n)_m(\gamma -\psi +n)_m}\prod_{j=1}^{k^*} s_j!,\\
\end{equation}
(ON): $s < m$ balls in $k^*$ {\it new} boxes in configuration $(s_1,\dots,s_{k^*})$ and the remaining $m-s$ balls in the $k$ {\it old} boxes in configuration $(m_1,\dots, m_k)$ for $\sum_{j=1}^{k} m_j= m-s$, $1 \leq s \leq m$, $\sum_{j=1}^{k^*} s_j=s$, $m_j \geq 0$, $s_j \geq 1$ with probability
\begin{equation}
\label{oldenew}
p_{{s, m}}({\bf n}):=\frac{(\gamma +n -k)_{m-k^*} (k- \psi)_{k^*} (k-\gamma +\psi)_{k^*}}{(\psi +n)_m(\gamma -\psi +n)_m}\prod_{j=1}^k (n_j+1)_{m_j\uparrow}\prod_{j=1}^{k^*}s_j.\\\\
\end{equation}	
Corresponding (O) and (N) one-step allocation rules under the new $(\gamma, \psi)$ parametrization arise respectively from (\ref{gibbsallold}) for $m=1$, $m_j=1$ and $m_l=0$ for $l \neq j$, and from (\ref{gibbsallnew}) for $m=1$, $k^*=1$ and $s_1=1$, and are given by
$$
\mbox{(O)}: p_{j}({\bf n}): =\frac{(n-k+\gamma)(n_j +1)}{n^2 +n\gamma +\psi(\gamma -\psi)} \mbox{ \hspace{0.2cm} for $j=1,\dots, k$  \hspace{0.0cm} and  \hspace{0.0cm} } \mbox{(N)}: p_0({\bf n}):=\frac{k^2 -k\gamma +\psi(\gamma -\psi)}{n^2 +n\gamma +\psi(\gamma -\psi)}.
$$\\
{\bf Remark 3.} Notice that (\ref{oldenew}), which is obtained specializing (19) in Cerquetti (2008),
provides, once the notation is made consistent, the explicit form for $p_{\bf b}(n_1, \dots, n_{\xi})$ in Sect. 7 of Gnedin (2010a) for the two-parameter $(\gamma, \psi)$ model. This kind of conditional distributions play a significant role in a Bayesian nonparametric approach to the treatment of species sampling problems under Gibbs priors, (see e.g. Favaro et al. 2009, Lijoi et al. 2007, 2008). Here we don't deal with this kind of applications. Some results in this perspective for the one parameter $(\gamma)$ Gnedin-Fisher model are in Cerquetti (2010).  

\section {Mixture representation and the number of occupied boxes}
The fundamental result in Gnedin and Pitman (2006, cfr. Th. 12) establishes that the EPPF of each Gibbs partition of genius $\alpha \in (-\infty, 1)$ corresponds to a mixture of extreme partitions probability function, which differ for $\alpha \in (-\infty, 0)$, $\alpha =0$ and $\alpha \in (0,1)$. Gnedin (2010a) provides the mixing law (\ref{mix}) over Poisson-Dirichlet $(-1, \xi)$ (i.e. $\alpha=-1$) extreme partitions  that corresponds to the parametrization (\ref{gnedintwo}). Here we derive the mixing law  for the reparametrization introduced in Theorem 1. as the limit distribution of the number of blocks following the approach in Gnedin (2010a). Additionally, by an application of Bayes theorem, we provide a direct proof of the weights in model (\ref{miogned}) actually arising by mixing over $\xi$ the extreme $PD(-1, \xi)$ weights.\\\\ 
First notice that both the {\it prior} (\ref{priorone}) and the {\it posterior} (\ref{postgamma}) for the number of blocks of the one-parameter $(\gamma)$-Gnedin-Fisher model  may be rewritten respectively as
\begin{equation}
\label{priGW}
\mathbb{P}_\gamma(\Xi=\xi)=\frac{\gamma(1-\gamma)_{\xi-1}}{\xi!}=\frac{(1)_{\xi -1}}{\Gamma(\xi)} \frac{(1-\gamma)_{\xi-1} (\gamma)_1}{(1)_{\xi}}
\end{equation}
and
\begin{equation}
\label{postGW}
\mathbb{P}_\gamma(\Xi=\xi|K_n=k)= \frac{(k)_{\xi-1}}{\Gamma(\xi)} \frac{(k-\gamma)_{\xi-1} (n+\gamma -k)_k}{(n)_{k+\xi-1}},
\end{equation}
for $\xi=1,2,\dots$, which shows that both belong to the class of  {\it shifted} univariate {\it generalized Waring distributions}, (also known as inverse Markov-Polya). This is a family of distributions on $\mathbb{N} \cup 0$ (Irwin, 1975; Xekalaki, 1983; see also Johnson et al. 2005), whose probability mass function is given by
$$
\mathbb{P}(N=i)=\frac{(\rho)_\eta}{i!} \frac{(a)_i (\eta)_i}{(a +\rho)_{\eta+i}} 
$$
for $i=0,1,2, \dots,$ for parameter $a, \eta, \rho$ positive reals, which arises by $Beta(\rho, a)$ mixture of a Negative Binomial distribution  $(\eta, p)$. Hence equations (\ref{priGW}) and (\ref{postGW}) correspond respectively to  $NB(1, p)$ and $Be(\gamma, 1-\gamma)$ and $NB(k, p)$ and $Be(n +\gamma -k, k-\gamma )$.
The probability generating function is, except for the constant, the Gaussian hypergeometric function
$$
_2F_1(a,\eta; a+\eta+\rho;z)=\sum_{i=0}^{\infty} \frac{(a)_i(\eta)_i}{(a +\eta +\rho)_i}\frac{z^i}{i!}
$$
which implies the generalized Waring distribution is overdispersed and characterized by heavy tail effect. Moreover $E(X^k)< \infty$ if and only if $\rho >k$. \\\\
The following result shows the reparametrization introduced in Theorem 1. yields even for the two-parameter $(\gamma, \psi)$ Gnedin-Fisher model a representation in terms of shifted generalized Waring mixture of Fisher $(-1, \xi)$ models.\\\\
{\bf Theorem 4.} {\it The EPPF in} (\ref{miogned}) {\it arises by mixing the family of $PD(-1, \xi)$ partition models 
\begin{equation}
\label{alfa1}
p_{\xi, -1}(n_1, \dots, n_k)=\frac{(\xi-1)_{k-1\uparrow -1}}{(\xi+1)_{n-1}}\prod_{j=1}^{k} {n_{j}},
\end{equation} 
over $\xi$, with a shifted generalized Waring distribution of parameters $a= 1-\gamma +\psi$, $\eta= 1-\psi$ and $\rho=\gamma$}, 
\begin{equation}
\label{mixlaw}
\mathbb{P}_{\gamma, \psi}(\Xi=\xi)= \frac{(1 -\psi)_{\xi-1}(1-\gamma +\psi)_{\xi-1}(\gamma)_{1-\psi}}{\Gamma(\xi) (1+\psi)_{\xi -\psi}},
\end{equation}
for $\psi \in [0,1)$ and $\gamma \in (0, \psi+1)$. \\\\
{\it Proof.}
By the general formula (\ref{nblocks}) for the law of the number of blocks for Gibbs partitions of genius $\alpha$, and exploting the definition of Lah numbers, the analogous of (\ref{number}) for the two-parameter model is given by
$$
\mathbb{P}_{\gamma, \psi}(K_n=k)= {n-1 \choose k-1} \frac{n!}{k!} \frac{(\gamma)_{n-k} (1- \psi)_{k-1} (1 -\gamma +\psi)_{k-1}}{(1+\psi)_{n-1}(1 +\gamma -\psi)_{n-1}}.
 $$
Rewriting in terms of Gamma functions yields
$$
\mathbb{P}_{\gamma, \psi}(K_n=k)=\frac{\Gamma(1+\psi) \Gamma(1+\gamma -\psi) (1-\psi)_{k-1} (1- \gamma +\psi)_{k-1}}{\Gamma(k+1) \Gamma(k) \Gamma(\gamma) } \frac{\Gamma(n+1) \Gamma(n) \Gamma(\gamma +n -k)}{\Gamma(n -k +1) \Gamma(\psi +n) \Gamma(\gamma -\psi +n)}
$$
and, by Stirling approximations, for $n \rightarrow \infty$ reduces to
\begin{equation}
\label{mixing2}
\mathbb{P}_{\gamma, \psi}(\Xi=\xi)= \frac{(1 -\psi)_{\xi-1}(1-\gamma +\psi)_{\xi-1}(\gamma)_{1-\psi}}{\Gamma(\xi) (1+\psi)_{ \xi-\psi }} 
\end{equation}
which provides the analogous of Eq. (5) in Gnedin (2010a) for the new parametrization. \\\\  
As $\xi \rightarrow \infty$ the power-like decay of the masses in (\ref{mixlaw}) (cfr. Gnedin, 2010a, Sect. 3) is rewritten as 
$$
\mathbb{P}_{\gamma, \psi} (\Xi = \xi ) \sim \frac{c}{\xi^{1+\gamma}}  \mbox{     with    } c=\frac{\Gamma(1+ \gamma -\psi) \Gamma(1 +\psi)}{\Gamma(1 -\psi) \Gamma(1- \gamma +\psi) \Gamma(\gamma)}.
$$\\
To show that the weights in (\ref{miogned}) actually arise by mixing the weights of the extreme Poisson-Dirichlet $(-1, \xi)$ partitions 
over $\xi$ with (\ref{mixing2}) we apply Bayes theorem. The {\it posterior} distribution for $\Xi$ for the reparametrized model may be obtained by the general form for Gibbs models of the posterior of the number of new blocks arising in a new sample of dimension $m$  (cfr. Cerquetti (2008, eq. (32) , see also Lijoi et al. 2007, eq. (4)) which expressed in terms of non-central generalized Stirling numbers is given by
\begin{equation}
\label{genpost}
\mathbb{P}(K^*_m=k^*| K_n=k)= \frac{V_{n+m, k+k^*}}{{V_{n,k}}}S_{m, k^*}^{-1,-\alpha, -(n -\alpha k)}.
\end{equation}
For $\alpha=-1$, inserting the specific weights in (\ref{miogned}), and exploiting the definition of non-central Lah numbers  $S_{n,k}^{-1, 1, r}= \frac{n!}{k!} {n -r-1 \choose n-k}$, (\ref{genpost}) yields
\begin{equation}
\label{postfinite}
\mathbb{P}_{\gamma, \psi} (K_m=k^*|K_n=k) = {m \choose k^*} (\gamma +n -k)_{m-k^*} (n+k+k^*)_{m-k^*} \frac{(k - \psi)_{k^*}(k -\gamma +\psi)_{k^*}}{(n +\psi)_{m} (n +\gamma -\psi)_m},
\end{equation}
and for $m \rightarrow \infty$, by standard Stirling approximations, the  {\it posterior} for the number of blocks of the two-parameter $(\gamma, \psi)$ model results
\begin{equation}
\label{posttwo}
\mathbb{P}_{\psi, \gamma}(\Xi=\xi |K_n=k)=\frac{(k-\psi)_{\xi-1} (k -\gamma +\psi)_{\xi-1}(n+\gamma -k)_{k-\psi}}{\Gamma(\xi) (n+\psi)_{k-\psi+\xi-1}},
\end{equation}
which is still in the class of shifted univariate Waring distributions for parameters $a= k-\gamma +\psi$, $\eta= k-\psi$ and $\rho= n+\gamma -k$, for $\gamma < k < n +\gamma$. Now, by Bayes theorem,
$$
\mathbb{P}_{\psi,\gamma}(\Xi=\xi) V_{n,k}^{-1, \xi}=\mathbb{P}_{\psi, \gamma}(\Xi=\xi |K_n=k) V_{n,k}^{\psi, \gamma},
$$
therefore, exploiting (\ref{posttwo})
$$ 
V_{n,k}^{\psi, \gamma}=\frac{(1 -\psi)_{\xi-1}(1-\gamma +\psi)_{\xi-1}(\gamma)_{1-\psi}}{\Gamma(\xi) (1+\psi)_{1-\psi +\xi -1}} \frac{(\xi-1)_{k-1\uparrow -1}}{(\xi+1)_{n-1}}\frac{\Gamma(\xi) (n+\psi)_{k-\psi+\xi-1}}{(k-\psi)_{\xi-1} (k -\gamma +\psi)_{\xi-1}(n+\gamma -k)_{k-\psi}}=
$$
and with the substitution $\xi -k=y$, 
$$ 
=\frac{(1 -\psi)_{ y+k-1}(1-\gamma +\psi)_{y+k-1}(\gamma)_{1-\psi}}{\Gamma(y+k) (1+\psi)_{y +k -\psi }} \frac{(y+k-1)_{k-1\uparrow -1}}{(y+k+1)_{n-1}}\frac{\Gamma(y+k) (n+\psi)_{k-\psi+y+k-1}}{(k-\psi)_{y+k-1} (k -\gamma +\psi)_{y+k-1}(n+\gamma -k)_{k-\psi}}.
$$
Then, by the multiplicative property of rising factorials $(x)_{a+b}=(x)_a (x+a)_b$, the last expression easily simplifies to 
$$
V_{n,k}^{\psi, \gamma}=\frac{(\gamma)_{n-k}(1 -\gamma +\psi)_{k-1} (1 -\psi)_{k-1}}{(1 +\psi)_{n-1} (1 +\gamma -\psi)_{n-1}},
$$
and the proof is complete. \hspace{11cm}$\square$\\\\
We provide an additional result for the two-parameter $(\gamma, \psi)$ Gnedin-Fisher model by exploiting the mixture representation introduced in Theorem 4. to obtain the {\it structural distribution}, the law of the frequency of Box 1,
$$
\tilde{P}_{1}:=  \lim_{n \to \infty } \frac{\#(B_1 \cap [n])}{n}.
$$\\
{\bf Proposition 5.} {\it The frequency $\tilde{P_1}$ of box $B_1$ of the $(\gamma, \psi)$ Gnedin-Fisher model has distribution
\begin{equation}
\label{lawP_1}
\mathbb{P}_{\gamma, \psi}(\tilde{P}_1 \in dy)= \gamma_{1-\psi} \Gamma(1+\psi) \left[\delta_1(dy)+(1-\gamma +\psi)(1-\psi) _2F_1 (2-\psi, 2-\gamma+\psi, 2; 1-y)
\right]dy
\end{equation}
for $_2F_1(a, b, c; x)$ the Gaussian hypergeometric function.}\\\\
{\it Proof:} By the mixture representation of Theorem 4.
$$
\mathbb{P}_{\gamma, \psi}(\tilde{P}_1 \in dy)= \sum_{\xi=1}^\infty \mathbb{} \mathbb{P}_{\gamma, \psi}(\Xi=\xi) \mathbb{P}(\tilde{P}_{\xi, 1} \in dy).
$$
By the theory of the symmetric Dirichlet model, it is known that $\tilde{P}_{\xi, i}\stackrel{d}=Beta(2, \xi-i)$,
therefore, since $Be(2, 0)=\delta_1(dy)$ 
$$
\mathbb{P}_{\gamma, \psi}(\tilde{P}_1 \in dy)=(\gamma)_{1 -\psi} \Gamma(1+\psi)\delta_1 (dy) +\sum_{\xi=2}^\infty \frac{(1 -\psi)_{\xi-1}(1-\gamma +\psi)_{\xi-1}(\gamma)_{1-\psi}}{\Gamma(\xi) (1+\psi)_{\xi -\psi}} \frac{\Gamma(\xi +1)}{\Gamma(\xi -1)}y (1-y)^{\xi -2}.
$$
By the change of variable $\xi-2=z$
$$
\mathbb{P}_{\gamma, \psi}(\tilde{P}_1 \in dy)=(\gamma)_{1 -\psi} \Gamma(1+\psi)\delta_1(dy) +\sum_{z=0}^\infty
\frac{(1-\psi)_{z+1}(1-\gamma +\psi)_{z+1} (\gamma)_{1-\psi}\Gamma(z+3)}{\Gamma(z+1)\Gamma(z+2)(1+\psi)_{z+2-\psi}} y(1-y)^{z}=
$$
and by standard combinatorial calculus
$$
=(\gamma)_{1 -\psi} \Gamma(1+\psi)\left[\delta_1(dy) + (1-\psi)(1 -\gamma +\psi)\sum_{z=0}^\infty
\frac{(2-\psi)_{z}(2-\gamma +\psi)_{z} }{\Gamma(z+1)\Gamma(z+2)} y(1-y)^{z}\right]
$$
and the result follows.\hspace{11.7cm}$\square$\\\\
{\bf Remark 6.} For $\psi=0$  (\ref{lawP_1}) yields the result in Gnedin (2010a, Sect. 6) for the distribution of the frequency $\tilde{P}_1$ of box $B_1$ for the one parameter $(\gamma)$ Gnedin-Fisher model. In fact
$$
_2F_1(2, 2-\gamma, 2; 1-y)= \sum_{z=0}^\infty \frac{(2-\gamma)_z}{\Gamma(z+1) (2)_z}(1-y)^z, 
$$
and multiplying and dividing by $y^{1-\gamma}$ and exploiting the probability mass function of the Negative Binomial $(2-\gamma, y)$ yields
$$
\sum_{z=0}^\infty \frac{(2-\gamma)_z}{\Gamma(z+1) (2)_z}y(1-y)^z = y^{\gamma -1},
$$
hence for $y \in (0,1]$
$$
\mathbb{P}_\gamma(\tilde{P}_1 \in dy)= \gamma \delta_1(dy)+ \gamma(1 -\gamma) y^{\gamma-1}dy.
$$



\section*{References}
\newcommand{\bibu}{\item \hskip-1.0cm}
\begin{list}{\ }{\setlength\leftmargin{1.0cm}}

\bibu \textsc{Cerquetti, A.} (2008) Generalized chinese restaurant construction of exchangeable Gibbs partitions and related results. {arXiv:0805.3853v1 [math.PR]}

\bibu \textsc{Cerquetti, A.} (2010) Bayesian nonparametric analysis of a species sampling model with finitely many types. {arXiv:1001.0245v1 [math.PR]}

\bibu  \textsc {Charalambides, C. A.} (2005) {\it Combinatorial Methods in Discrete Distributions}. Wiley, Hoboken NJ.

\bibu \textsc{Favaro, S., Lijoi, A., Mena, R. and Pr\"unster, I.} (2009) Bayesian non-parametric inference for species variety with a two-parameter Poisson-Dirichlet process prior. {\it JRSS-B}, 71, 993-1008.

\bibu \textsc{Fisher, R.A., Corbet, A.S. and Williams, C.B.} (1943) The relation between the number of species and the number of individuals in a random sample of an animal population. {\it J. Animal. Ecol.}, 12, 42--58.

\bibu \textsc{Gnedin, A.} (2010a) A species sampling model with finitely many types. {\it Elect. Comm. Probab.}, 15, 79--88.

\bibu \textsc{Gnedin, A.} (2010b) Boundaries from inhomogeneous Bernoulli trials. {arXiv:0909.4933 [math.PR]} 

\bibu \textsc{Gnedin, A. and Pitman, J. } (2006) {Exchangeable Gibbs partitions  and Stirling triangles.} {\it J. Math. Sci}, 138, 3, 5674--5685. 

\bibu \textsc{Hansen, B. and Pitman, J.} (2000) Prediction rules for exchangeable sequences related to species sampling. {\it Stat. \& Probab. Letters}, 46, 251--256.

\bibu \textsc{Irwin, J.O.} (1975) The Generalized Waring distribution. Part I. {\it J. R. Statist. Soc. A}, 138, p.18

\bibu \textsc{Johnson, N.L., Kotz, S., Kemp, A.W.}(2005) {\it Univariate Discrete Distributions}, 3rd Ed. Wiley, New YOrk, 2005.

\bibu \textsc{Lijoi, A., Mena, R.H. and Pr\"unster, I.} (2007) Bayesian nonparametric estimation of the probability of discovering new species.  {\it Biometrika}, 94, 769--786.
 
\bibu \textsc{Lijoi, A., Pr\"unster, I. and Walker, S.G.} (2008) Bayesian nonparametric estimator derived from conditional Gibbs structures. {\it Annals of Applied Probability}, 18, 1519--1547.

\bibu \textsc{Pitman, J.} (2006) {\it Combinatorial Stochastic Processes}. Ecole d'Et\'e de Probabilit\'e de Saint-Flour XXXII - 2002. Lecture Notes in Mathematics N. 1875, Springer.

\bibu \textsc{Pitman, J. and Yor, M.} (1997) The two-parameter Poisson-Dirichlet distribution derived from a stable subordinator. {\it Ann. Probab.}, 25:855--900.

\bibu \textsc{Xekalaki, E.} (1983) Infinite divisibility, completeness and regression properties of the univariate generalized Waring distribution. {\it Ann. Inst. Statist. Math.}, 35, 279-289

\end{list}
\end{document}